\documentclass[twoside]{univzil}
\usepackage{amsmath}
\usepackage{amssymb,amsfonts}
\usepackage{graphicx}                 %% for include pictures *.eps
\usepackage{psfrag}                   %% for using substitution \psfrag
\begin{document}
\renewcommand\volumeZA{20/2005}          %% volume of journal

%============================================================================%
%                              head of article                               %
%============================================================================%

%%%%%%%%%%%%%%%%%%%%%%%%%%%%%%%%%%%%%%%%%%%%%%%%%%%%%%%%%%%%%%%%%%%%%%%% title
\title[Line transversals and homothety transforms]{Line transversals \\
for homothetical systems of polygons in ${\mathbb R}^2$}

%%%%%%%%%%%%%%%%%%%%%%%%%%%%%%%%%%%%%%%%%%%%%%%%%%%%%%%%%%%%%%%%%%%%%% authors
%% First, second, third, ...  author (name, address, email)
%% If there are more authors then all authors contain the same items
%% \author{Name}
%% \address{Street XY, Town, State}  \curraddress{current adress}
%% \email{mail@aaa.bbb.com}          \urladdress{http://...}

\author{M. Kauki\v c}
\address{Dept. of Mathematical Methods, FRI \v ZU, Ve\v lk\'y Diel, 
O1026 \v Zilina, Slovakia}
\email{mike@frcatel.fri.utc.sk}

%%%%%%%%%%%%% (optional) If any thanks for the financial supports, grants, ...
%%%%%%%%%%%%%%%%%%%%%%%%%%%%%%%% \thanks{The first  author was supported ... }
\thanks{This research was supported by Slovak Scientific
Grant Agency under grant No. 1/0490/03}

%%%%%%%%%%%%%%%%%%%%%%%%%%%%%%%%%%%%%%%%%%%%%%%%%%%%%%%%%% (optional) Keywords

%%%%%%%%%%%%%%%%%%%%%%%%%%%%%%%%%%%%%%%%% (obligatory) AMS Classification 2000
%% The 1-st classification is obligatory, the 2-nd classification is optional
%% \subjclass{primary}{secondary}      f.e. \subjclass{35R35, 49N50}{}
%%%%%%%%%%%%%%%%%%%%%%%%%%%%%%%%%%%%%%%%%%% \subjclass{35R35, 49M15}{49N50}
{}
%%%%%%%%%%%%%%%%%%%%%%%%%%%%%%%%%%%%%%%%%%%%%%%%%%%%%%%%% (optional) Abstract
\begin{abstract}
The problem considered in this article is: `for given finite system 
of convex polygons in the plane which have no transversal, 
find such homothety transformations of polygons (having fixed centres 
inside given polygons) with minimal similarity ratio $c>1$ 
that the transformed system has a
transversal'. We prove that in this `minimal configuration', 
we can always find three polygons (two of them lying in distinct halfplanes
determined by the transversal), for which the transversal is also the
tangent line.
\end{abstract}

%%%%%%%%%%%%%%%%%%%%%%%%%%%%% private macros (f.e. the different environments)
\newtheorem{theorem}{Theorem}[section]
\newtheorem{corollary}[theorem]{Corollary}
\newtheorem{lemma}[theorem]{Lemma}
\newtheorem{exmple}[theorem]{Example}
\newtheorem{defn}[theorem]{Definition}
\newtheorem{proposition}[theorem]{Proposition}
\newtheorem{conjecture}[theorem]{Conjecture}
\newtheorem{rmrk}[theorem]{Remark}
            %%% for no-italic, numbered environments, use:
\newenvironment{definition}{\begin{defn}\normalfont}{\end{defn}}
\newenvironment{remark}{\begin{rmrk}\normalfont}{\end{rmrk}}
\newenvironment{example}{\begin{exmple}\normalfont}{\end{exmple}}
            %%% for unnumbered environments, use f.e.
\newtheorem*{remarque}{Remark}
\newcommand\refrm[1]{{\rm (\ref{#1})}}
\newcommand\bbreak{\allowdisplaybreaks}
\newcommand\dd{\mathop{\rm d\!}\nolimits}
\newcommand\sgn{\mathop{\rm sgn}\nolimits}

%%%%%%%%%%%%%%%%%%%%%%%%%%%%%%%%%%%%%%%%%%%%%%%%%%%%%%%%%%%%%%%%%%%% maketitle
\maketitle

%============================================================================%
%                             the main article                               %
%============================================================================%

%%%%%%%%%%%%%%%%%%%%%%%%%%%%%%%%%%%%%%%%%%%%%%%%%%%%%%%%%%%%%%%%%%%%%%%%%%%%%%
% homothety with center O and ratio k - terminologia

\section{Introduction, notation}

\begin{definition}
The {\bf line transversal} to a system of convex sets in 
$\mathbb{R}^d,\, d\ge 2$ is 
the line, which intersects every member of the system.
\end{definition}
For brief
overview of recent results in geometric transversal theory,
see \cite{wenger}. In this article we restrict our attention 
to the convex polygons in the plane. 

\begin{definition}
We will call the finite system $\mathcal{S}=\{P_1,P_2,\dots,P_n\}$ 
of convex polygons in $\mathbb{R}^2$ the {\bf initial configuration} 
if no transversal of $\mathcal{S}$ exists.
\end{definition}

\begin{figure}[h]
\centering
\psfrag{t}{$t$}
\includegraphics[width=8.5cm,height=4.8cm]{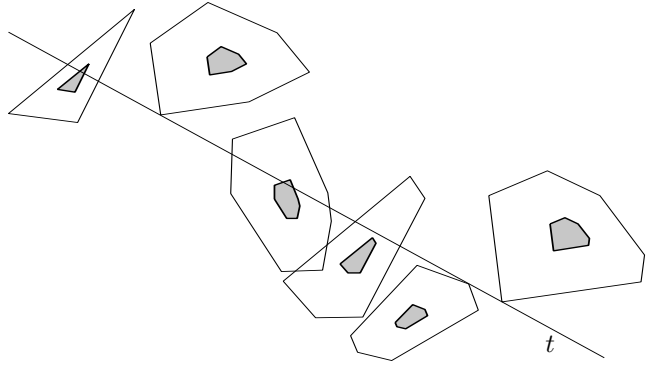}
\caption{\bf Initial configuration and its minimal configuration}
\label{pic1}
\end{figure}

For example, the smaller, filled polygons on Fig. \ref{pic1} form an initial
configuration. 

In this paper we will use the notation $\mathcal{H}_{S,c}$ for the homothety
transformation with center $S$ and similarity ratio $c$ given by:
$$
\overrightarrow{\mathcal{H}_{S,c}(X)-S} = c\, \overrightarrow{(X-S)},\, 
\mbox{ or }\, 
\mathcal{H}_{S,c}(X) = S + c\, \overrightarrow{(X-S)},
\ X \in \mathbb{R}^2.
$$
Let $P$ be the convex polygon and $S\in P$ an interior point of $P$.
Let us denote by 
$$\mathcal{H}_{S,c}(P)=\{\mathcal{H}_{S,c}(X),\, X \in P \}$$
the image of polygon $P$ in transformation $\mathcal{H}_{S,c}$.
Throughout this paper we will fix the center $S$ of homothety transformation 
applied to polygon $P$ to be the centroid of $P$, i.e. the point $T=(x_T,y_T)$,
with coordinates $x_T$, resp. $y_T$ determined by the arithmetic means of 
$x\/$--coordinates (resp. $y\/$--coordinates) of vertices of $P$. 
In this case we will use the simple notation $\mathcal{H}_c(P)$ 
for the image of $P$.  

We use the notation $\mathcal{S}_c$ for the image of finite system
$\mathcal S=\{P_1, P_2,\dots ,P_n\}$ of polygons,
obtained by applying the transformations $\mathcal{H}_c$ to every polygon
from $\mathcal{S}$. For $c>1$ (resp. $c<1$) 
the resulting image $\mathcal{S}_c$ can be described as the configuration 
consisting of  polygons of $\mathcal{S}$ expanded (resp. shrinked) 
around their centroids as can be seen on  Fig. \ref{pic1} (bigger polygons
drawn by thin lines).

It is clear that every finite system $\mathcal S$ of convex polygons 
in the plane can be transformed to an initial configuration $\mathcal S_c$ 
by taking the similarity ratio $c$ sufficiently small, 
except for the degenerate case, 
where all the centroids of polygons of $\mathcal S$ are collinear. 

For given initial configuration $\mathcal S$ there exists a real number
$c_m > 1$ such that
for every positive $c<c_ m$ the configuration $S_c$ has no line traversal  
and the configuration $S_{c_m}$ has a line traversal.

\begin{definition}
We will call the number $c_m$ the {\bf minimal expansion ratio}
and the configuration ${\mathcal S}_{c_m}$ the {\bf minimal configuration} for the
given initial configuration $\mathcal S$.
\end{definition}

Example of minimal configuration (big non-filled polygons) for initial
configuration (given by small filled polygons) can be seen on
Fig. \ref{pic1}.

\section{Main result}

We can see from Fig. \ref{pic1} that, in minimal configuration, there are
three polygons, for which the (unique) transversal $t$ is also the tangent line.
This is the general property of minimal configuration, 
as we state in the next theorem. 

\begin{theorem}
\label{mainthm}
Let $\mathcal S$ be the initial configuration. The configuration 
${\mathcal S}_{c_m}$ with $c_m>1$ is the minimal configuration 
for $\mathcal S$ if and only if there are three
polygons $P_i, P_j, P_k \in {\mathcal S}_{c_m}$, which  intersect 
the transversal $t$ to ${\mathcal S}_{c_m}$ either in one point (vertex) 
or have one common side with $t$. Polygons $P_i, P_j, P_k$ cannot lie 
in the same halfplane determined by transversal $t$. The transversal $t$
for minimal configuration is unique. 
\end{theorem}

Before the proof of theorem \ref{mainthm} let us make some remarks
about concepts, involved in the proof.

\begin{definition}
Let us have the line $q$ and the convex polygon $P$. Then there exists 
a nonnegative real number $c$, such that the line $q$ is the tangent line 
of transformed polygon $\mathcal{H}_c(P)$. We will call such $c$
the {\bf correcting factor of polygon $P$ with respect to line} $q$. 
If the centroid of $P$ lies on the line $q$, we put $c=0$.
\end{definition}

\begin{figure}[h]
\centering
\psfrag{P}{$P$}
\psfrag{L1}{$q_1$}
\psfrag{L2}{$q_2$}
\psfrag{S}{$S$}
\psfrag{Pl}{$P_r$}
\psfrag{Pu}{$P_u$}
\psfrag{Ll}{$q_r$}
\psfrag{Lu}{$q_u$}
\includegraphics[width=5.5cm]{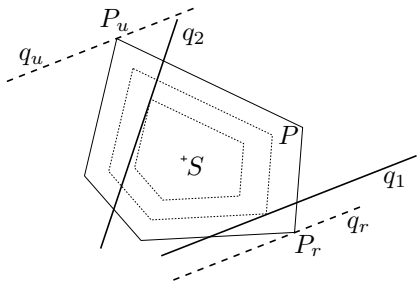}
\caption{\bf Shrinking factor}
\label{pic2}
\end{figure}

We will use the above definition mainly in the case, when the line $q$
has nonempty intersection with interior of the polygon $P$. In such case,
we will speak also about {\bf shrinking factor}, see Fig. \ref{pic2}.
For the line $q_1$ the shrinking factor is $c_1=3/4$ and for $q_2$
the shrinking factor is $c_2=1/2$.

It is not difficult to express the correcting factor explicitly.
For given line $q$ and polygon $P$ we can construct the `minimal strip'
containing $P$ with boundary lines, which are the translates of $q$
(see Fig \ref{pic2}, for the line $q_1$ the strip is bounded by the lines
$q_u, q_r$). Let us denote by $\varrho(A,p)$ the distance of point $A$
from the line $p$ and by $q_s$ the line, parallel with $q$ such that
$S \in q$ then the correcting factor can be expressed as:
 $$c=\left\{{\begin{array}{ll}
  \displaystyle{\frac{\varrho(S,q)}{\varrho(S,q_u)}}
       & \textrm{if $q,q_u$ lie in the same halfplane, determined by $q_s$} \\ 
  \displaystyle{\frac{\varrho(S,q)}{\varrho(S,q_r)}} 
       & \textrm{if $q,q_r$ lie in the same halfplane, determined by $q_s$.} 
 \end{array}
 }
 \right .$$
Thus, the following proposition holds. 
\begin{proposition}
Let us take $t=\varrho(S,q)$ as independent variable. 
Then the function $c=c(t),\, t\in \mathbb R$ is continuous,
piecewise linear, V-shaped function.
\label{prop1}
\end{proposition}

\begin{proof}[Proof of Theorem~{\upshape\ref{mainthm}}]
Let us consider the initial configuration $\mathcal S$.
We will show that the configuration ${\mathcal S}_c$, in which less than three 
polygons have the transversal as tangent line, cannot be minimal.

\noindent
{\bf Case 1.}
If, for some configuration ${\mathcal S}_c=\{P_1, P_2,\dots ,P_n\}$, 
there exists a transversal $t$,
which has nonempty intersections with interiors of all polygons  
from  ${\mathcal S}_c$, then all shrinking factors $c_1, c_2, \dots, c_n$
have values less than 1, and so has also the number
$$c_m = \max\{c_1, c_2, \dots, c_n \} <1 \, \mbox{ and we have }
\tilde{c} = c\, c_m < c.$$ 
Then, the configuration ${\mathcal S}_{\tilde{c}}$ has the same transversal
$t$, but smaller similarity ratio, so the configuration ${\mathcal S}_c$ 
cannot be minimal.

\noindent
{\bf Case 2.} Now let us consider the second case, 
when for some configuration 
${\mathcal S}_c$ there is exactly one polygon $P_k \in {\mathcal S}_c$ with 
no interior points belonging to given transversal $t$ of ${\mathcal S}_c$
(i.e. $t$ is the tangent line of $P_k$). Let us introduce the function
$\varphi(c)$ defined by formula:
$$\varphi(c)=c - \max\{c_1(c), c_2(c), \dots, 
c_{k-1}(c),c_{k+1}(c),\dots, c_n(c) \},$$
where $c_j(c), j \ne k$ are the correcting factors of polygons $P_j \ne P_k$,
when the transversal $t$ is translated so that the shrinking factor of polygon
$P_k$ will become $c$ (see Fig. \ref{pic3}). 
\begin{figure}[h]
\centering
\psfrag{t1}{$t (c=1)$}
\psfrag{tc}{$c < 1$}
\psfrag{tk}{$c_{krit}$}
\psfrag{S}{$S$}
\includegraphics[width=7cm]{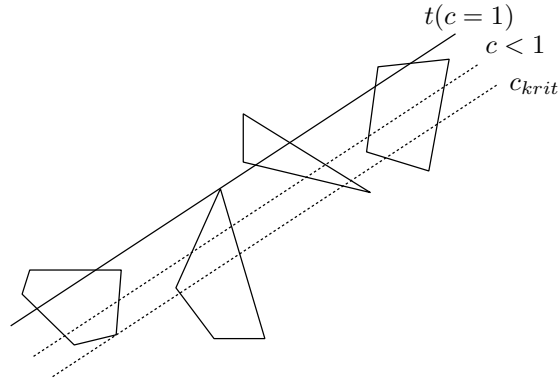}
\caption{\bf Case of one polygon touching transversal}
\label{pic3}
\end{figure}

For $c=c_k=1$ we have the untranslated transversal $t$. We will consider
only values of $c$ in interval $c_{krit} \le c \le 1$, where $c_{krit}$
is such number that for $c<c_{krit}$ the translated line $t$ remains no longer
the transversal of ${\mathcal S}_c$, so $c_j=1$ for some $j \ne k$ if
$c=c_{krit}$ 
(cf. Fig. \ref{pic3}). It is clear that 
$c_{krit} < 1$, because if it is not true, then some other polygon 
$P_i, i \ne k$ should touch the transversal $t$. 
It can be easily shown that the function $\varphi(c)$ is
continuous (it follows from proposition \ref{prop1}). Now, we have:
$$\varphi(1)=1-\max\{c_1, c_2, \dots, c_{k-1},c_{k+1},\dots, c_n \} >0,
\quad \varphi(c_{krit})=c_{krit}-1 < 0$$
it follows that for some $\tilde{c},\, c_{krit} < \tilde{c} < 1, \,
\varphi(\tilde{c})=0$ or 
$$\tilde{c} = c_k =
\max\{c_1(\tilde{c}), c_2(\tilde{c}), 
\dots, c_{k-1}(\tilde{c}),c_{k+1}(\tilde{c}),\dots, c_n(\tilde{c}) \}.$$    
For such $\tilde{c}$ the translated transversal is the also transversal
for configuration $\mathcal{S}_{c\,\tilde{c}}$, so the original configuration
$\mathcal{S}_c$ cannot be minimal.

\noindent
{\bf Case 3.} Finally, lets us consider some configuration 
${\mathcal S}_c$ in which there are exactly two polygons 
$P_i, P_j \in {\mathcal S}_c$ with 
no interior points belonging to given transversal $t$ of ${\mathcal S}_c$.
It is sufficient to take into account only such configuration, where
$P_i, P_j$ lie in distinct halfplanes, determined by the transversal $t$,
because if they are situated in the same halfplane, we can proceed as in
preceding case.
Now, consider the function:
$$\varphi(c)=c - \max\{c_k(c); \, k=1,2,\dots,n,\, k\ne i, k\ne j \},$$
where $c_k(c)$ are the correcting factors of polygons $P_k$,
when the transversal $t$ is moved so that the shrinking factor of both polygons
$P_i, P_j$ will become $c$ (Fig. \ref{pic4}). 
\begin{figure}[h]
\centering
\psfrag{t}{$t (c=1)$}
\psfrag{c1}{$c < 1$}
\psfrag{ck}{$c_{krit}$}
\psfrag{Pi}{$P_i$}
\psfrag{Pj}{$P_j$}
\includegraphics[width=7cm]{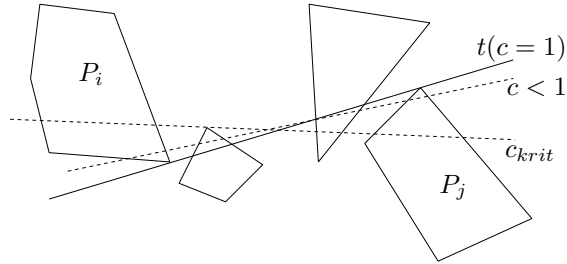}
\caption{\bf Case of two polygon touching transversal}
\label{pic4}
\end{figure}
 
Analogically as in case 2. we will consider only values of $c$ 
lying in interval $c_{krit} \le c \le 1$, where $c_{krit}$
is such number that for $c<c_{krit}$ the moved line $t$ remains no longer
the transversal of ${\mathcal S}_c$, so $c_l(c_{krit})=1$ for some
$l \ne i, l\ne j$.

If the function $\varphi(c)$ is continuous, we can proceed exactly as in 
case 2., showing that moved transversal $t$ is also the transversal for some
configuration ${\mathcal S}_{\hat c}$ with $\hat c < c$; thus, ${\mathcal S}_c$
cannot be the minimal configuration.  

It is sufficient to ascertain the continuity of individual coefficients
$c_k(c)$. Let us look at Fig. \ref{pic5}.
\begin{figure}[h]
\centering
\psfrag{t1}{$t_1$}
\psfrag{t2}{$t_2$}
\psfrag{Pk}{$P_k$}
\includegraphics[width=7cm]{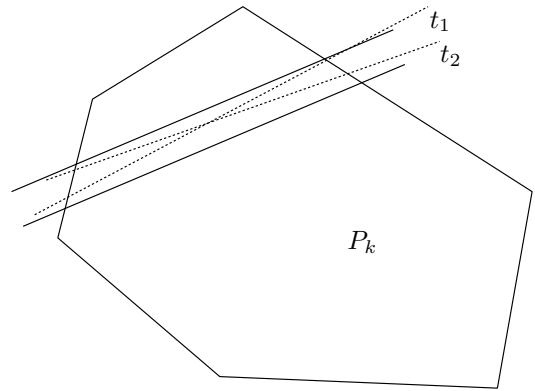}
\caption{\bf Two lines $t_1,t_2$ with slightly different values of $c$}
\label{pic5}
\end{figure}

When we consider two moved transversals $t_1, t_2$ for slightly different
values $c_1, c_2$ of $c$, it is clear, that the angle between them can be made
arbitrarily small by choosing the difference $|c_1 - c_2|$ to be sufficiently
small. Now, we can enclose both $t_1 \cap P_k$ and $t_2 \cap P_k$ in a strip, 
bounded by parallel lines. Because the width of the strip can be made also 
arbitrarily small by choosing $c_1, c_2$ to be very close, we have reduced
the problem of continuity of $c_k(c)$ to the analogical problem for translates,
which was analysed in case 2. 

So, considering above three cases, we prove that if ${\mathcal S}_c$ is minimal
configuration for $\mathcal S$ then (at least) three polygons from 
${\mathcal S}_c$ exist, which have the transversal $t$ to ${\mathcal S}_c$ 
as the tangent line. All three polygons cannot lie in the same halfplane,
determined by $t$, because such configuration can be easily reduced by method
used in case 2.

Conversely, if for some configuration ${\mathcal S}_c$ we can find three
polygons from ${\mathcal S}_c$, which have the transversal $t$ 
as common tangent and  these polygons do not lie in the same halfplane,
determined by $t$, then it is obvious (see Fig. \ref{three}; only three
mentioned polygons are shown, the location of other polygons from
${\mathcal S}_c$ does not matter) that 
${\mathcal S}_c$ is minimal configuration.
\begin{figure}[h]
\centering
\psfrag{t1}{$t_1$}
\psfrag{t2}{$t_2$}
\psfrag{Pk}{$P_k$}
\includegraphics[width=7cm]{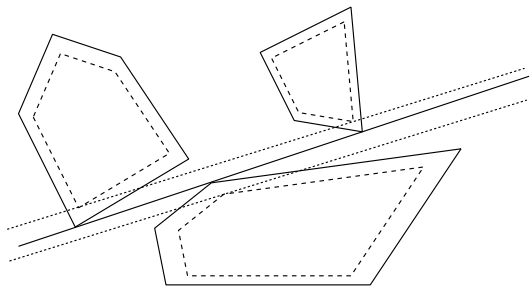}
\caption{\bf Configuration with three polygons touching transversal}
\label{three}
\end{figure}

Really, by any change of $c$ to smaller value $\tilde c$, we obtain
a strip of nonzero width, bounded by parallel lines, 
which separates the mentioned polygons,
so we cannot find the line (transversal), which has common points 
with all three polygons. From this fact it is also clear that the traversal
for minimal configuration is unique.
\end{proof}

%============================================================================%
%                                references                                  %
%============================================================================%

\end{document}